\documentclass{article}
\usepackage{amsfonts}
\title{\bf Some analytical results associated with extensions of the canonical Feller--Spitzer distribution}
\date{}  
 
\author{
\ R.B. Paris\footnote{Division of Computing and Mathematics, Abertay University, Dundee DD1 1HG, UK.\ \ 
E-Mail: r.paris@abertay.ac.uk}
\ \ and
Vladimir V. Vinogradov\footnote{Department of Mathematics, Ohio University, Athens, OH, USA\ \ 
E-mail: vinograd@ohio.edu}}

\begin{document}
\def\f#1#2{\mbox{${\textstyle \frac{#1}{#2}}$}}
\def\dfrac#1#2{\displaystyle{\frac{#1}{#2}}}
\def\boldal{\mbox{\boldmath $\alpha$}}
{\newcommand{\Sgoth}{S\;\!\!\!\!\!/}
\newcommand{\bee}{\begin{equation}}
\newcommand{\ee}{\end{equation}}
\newcommand{\la}{\lambda}
\newcommand{\ka}{\kappa}
\newcommand{\al}{\alpha}
\newcommand{\fr}{\frac{1}{2}}
\newcommand{\fs}{\f{1}{2}}
\newcommand{\g}{\Gamma}
\newcommand{\br}{\biggr}
\newcommand{\bl}{\biggl}
\newcommand{\ra}{\rightarrow}
\newcommand{\gl}{\raisebox{-.8ex}{\mbox{$\stackrel{\textstyle >}{<}$}}}
\newcommand{\gtwid}{\raisebox{-.8ex}{\mbox{$\stackrel{\textstyle >}{\sim}$}}}
\newcommand{\ltwid}{\raisebox{-.8ex}{\mbox{$\stackrel{\textstyle <}{\sim}$}}}
\renewcommand{\topfraction}{0.9}
\renewcommand{\bottomfraction}{0.9}
\renewcommand{\textfraction}{0.05}
\newcommand{\mcol}{\multicolumn}
\date{}
\maketitle
\pagestyle{myheadings}
\markboth{\hfill \sc R.B. Paris and V.V. Vinogradov  \hfill}
{\hfill \sc  Running title: Some results related to generalized Feller-Spitzer distributions \hfill}
\begin{abstract}

In this work, we establish new analytical results which are required for the derivation of subtle 
properties of the members of two classes of the generalized Feller--Spitzer distributions introduced in our paper \cite{PV6}.

\vspace{0.4cm}


\vspace{0.3cm}


\end{abstract}

\vspace{0.3cm}

\newtheorem{theorem}{Theorem}
\begin{center}
\noindent {\bf 1. \  Preliminaries}
\end{center}
\setcounter{section}{1}
\setcounter{equation}{0}
\renewcommand{\theequation}{\arabic{section}.\arabic{equation}}
This work contains several analytical results which are closely related to two classes of the generalized Feller--Spitzer distributions 
introduced in our recent paper \cite{PV6}. In particular, Theorem \ref{slozhnaya} provides a series expansion for the probability density function 
of a certain absolutely continuous distribution, 
whereas Lemmas \ref{Mera!Levy2} and \ref{Poincare}  present useful closed-form expression in terms of hypergeometric and digamma functions, and 
the Poincar\'{e} asymptotic series, 
respectively, for a certain {\it ``integrated tail"} which involves a specific modified Bessel function (\ref{ModBes1st}).
Finally, the new Bessel function inequalities (\ref{a1}) and (\ref{a2}) are employed in  \cite{PV6} for the derivation of  some properties of 
the generalized Feller--Spitzer distributions.

\bigskip 

Next, we introduce some relevant special functions that will be used in this work. 

\newtheorem{definition}{Definition}
\begin{definition}	\label{dsDef}
({\it Gauss hypergeometric function}, see, for example, \cite[(15.2.1)]{DLMF}). Given arbitrary complex $a$ and $b$,  $c \neq 0, -1, -2, \ldots$, 
and argument $z \in \mathbb{C}$ with $|z| < 1$,  the convergent series which emerges on the right-hand side of (\ref{DiStPGF})  
\begin{equation}	\label{DiStPGF}	
{}_{2}F_{1}(a, b; c; z) := \sum_{\ell = 0}^{\infty} \frac{(a)_{\ell} \cdot (b)_{\ell}}{(c)_{\ell}} \cdot 
\frac{z^{\ell}}{\ell!}
%
%
\end{equation}
generates the so-called {\it Gauss hypergeometric function}. Hereinafter, $(a)_k=\g(a+k)/\g(a)$ is the Pochhammer symbol.
\end{definition}

The infinite series which is present on the right-hand side of (\ref{DiStPGF}) can often be simplified. For instance, 
see  \cite[(4.1)]{PV4} for a few relevant special cases of function ${}_{2}F_{1}(1/2, 1; c; z)$.

\begin{definition}	\label{GENERdsDef}
({\it Generalized hypergeometric function} ${}_{3}F_{2}$, see \cite[(16.2.1) and case 16.2(iii)]{DLMF}). Given arbitrary complex $a$, $b$, $c$, $d  \neq 0, -1, -2, \ldots$, 
$e \neq 0, -1, -2, \ldots$, and argument $z \in \mathbb{C}$ with $|z| < 1$,  the function 
\begin{equation}	\label{generalizedDiStPGF}	
{}_{3}F_{2}\bl(\begin{array}{c} a, b, c\\d,e \end{array}; z\br) := 
\sum_{\ell = 0}^{\infty} \frac{(a)_{\ell} \cdot (b)_{\ell}  \cdot (c)_{\ell}}{(d)_{\ell} \cdot (e)_{\ell}} \cdot 
\frac{z^{\ell}}{\ell!}
\nonumber
\end{equation}
is called a {\it generalized hypergeometric function}  ${}_{3}F_{2}$.
\end{definition}

In addition, we refer to \cite[(16.2.1)]{DLMF} for a more general form ${}_{p}F_{q}$ of a generic member 
of the class of generalized hypergeometric functions.

The modified Bessel function of the first kind of order $\nu$ is defined by
\begin{equation} 	\label{ModBes1st}
I_\nu(t) := \sum_{k = 0}^{\infty} ~\frac{(t/2)^{2k + \nu}}{k! \cdot \g(k+1+ \nu)}\qquad (|t|<\infty).
\end{equation}
In the sequel, we will use the digamma function which is hereinafter denoted by 
\begin{equation} 	\label{diGAmma}
\psi(z) := \Gamma^{\prime}(z)/\Gamma(z)~~~\rm{where}~~~\Re(z) > 0
\end{equation}
with the Euler--Mascheroni constant $\gamma := - \psi(1)$.

\vspace{0.6cm}

\begin{center}
\noindent {\bf 2. \ A series expansion for a particular inverse Laplace transform}
\end{center}
\setcounter{section}{2}
\setcounter{equation}{0}
\renewcommand{\theequation}{\arabic{section}.\arabic{equation}}
Here, we concentrate on the consideration of the probability density function (introduced by (\ref{inverseLaTr})) of a certain non-negative random variable. 
Set
\bee		\label{funktsija-phi}
\phi(x) := x \cdot \sqrt{1-x^2}~~~~\rm{with~~~~{\it |x|}} \leq 1.
\ee
Observe that the trigonometric transformations of the function $\phi(x)$ defined by (\ref{funktsija-phi}) admit
the following representations in terms of Taylor series which are convergent for 
$|x| \leq 1$:
\[		\label{CosPhi}
\cos \phi(x)= \sum_{n\geq0} \frac{(-1)^n}{(2n)!}\,x^{2n} (1-x^2)^n,
\]
and
\[		\label{SinPhi}
\sin \phi(x)= \sum_{n\geq 0}\frac{(-1)^n}{(2n+1)!}\,x^{2n+1}(1-x^2)^{n+1/2}.
\]

For real $t > 0$, define 
\bee		\label{inverseLaTr}
{\mathcal I}(t) := \frac{e^{-t}}{2\pi i}\int_{c-\infty i}^{c+\infty i} \frac{e^{wt} e^{w\sqrt{w^2-1}-w^2+1}}{w+\sqrt{w^2-1}}\, \cdot dw.
\ee
Here, real $c > 1$ is fixed. 

It follows from Cauchy's theorem
that the inverse Laplace transform which emerges on the right-hand side of (\ref{inverseLaTr}) 
does not depend on a particular value of  real $c > 1$.
\begin{theorem}		\label{slozhnaya}
The following representation for the above function ${\mathcal I}(t)$ in terms of the 
difference of two convergent series holds: 
\[
{\mathcal I}(t) = \frac{2te^{-t}}{\sqrt{\pi}}\bl(\sum_{n\geq0}\frac{(-2)^n}{(2n)!} \sum_{r\geq0} b_r(n) \bl(\frac{1}{t}~ 
\cdot \frac{d}{dt}\br)^n [t^{-r-2} I_{2n+r+1}(t)]
\]
\bee		\label{okonchatPRED}
\hspace{1.6cm}- \sum_{n\geq0}\frac{(-2)^n}{(2n+1)!} \sum_{r\geq0} b_r(n) \bl(\frac{1}{t}~ \cdot \frac{d}{dt}\br)^{n+1} [t^{-r-1} I_{2n+r+2}(t)]\br), 
\ee
where $b_r(n) :=2^r \cdot \g(n+r+\f{3}{2})/r!$ for integer $r \geq 0$.
\end{theorem}

\noindent {\bf Proof of Theorem \ref{slozhnaya}.} We omit for now the factor $e^{-t}$ outside the integral 
in (\ref{inverseLaTr}), we evaluate the integral per se as follows:
\begin{eqnarray}		
\mathbf{I} &:=& \frac{1}{2\pi i}\int_{c-\infty i}^{c+\infty i} \frac{e^{wt} e^{w\sqrt{w^2-1}-w^2+1}}{w+\sqrt{w^2-1}}\,dw\nonumber\\
& = &\frac{1}{2\pi i}\int_{-1}^1 e^{xt-x^2+1}\bl\{\frac{e^{-ix\sqrt{1-x^2}}}{x-i\sqrt{1-x^2}}-\frac{e^{ix\sqrt{1-x^2}}}{x+i\sqrt{1-x^2}}\br\}dx\nonumber\\
& =&\frac{e}{\pi}\int_{-1}^1 e^{xt-x^2}\{\sqrt{1-x^2} \cos \phi(x)-x \sin \phi(x)\}dx.\label{samInteg}
\end{eqnarray}
In obtaining (\ref{samInteg}), the integration path $(c-\infty i,c+\infty i)$ has been collapsed onto the upper and lower sides of the branch cut in the $w$-plane between $[-1,1]$.

Subsequently, we can decompose the integral $\mathbf{I}$ given by (\ref{samInteg}) into the following difference:
\bee		\label{raznostINTEG}
\mathbf{I} = \mathbf{I}_1 - \mathbf{I}_2,
\ee
where
\begin{eqnarray}		
\mathbf{I}_1 &=& \frac{e}{\pi}\sum_{n\geq0}\frac{(-1)^n}{(2n)!}\int_{-1}^1 e^{xt-x^2}x^{2n}(1-x^2)^{n+1/2}dx\nonumber\\
& =&\frac{2e}{\pi}\sum_{n\geq0}\frac{(-1)^n}{(2n)!} \sum_{k\geq0}\frac{t^{2k}}{(2k)!} \int_0^1 e^{-x^2} x^{2n+2k}(1-x^2)^{n+1/2} dx \label{Iodin}
\end{eqnarray}
and
\begin{eqnarray}		
\mathbf{I}_2 &=& \frac{e}{\pi}\sum_{n\geq0}\frac{(-1)^n}{(2n+1)!}\int_{-1}^1 e^{xt-x^2}x^{2n+2}(1-x^2)^{n+1/2}dx\nonumber\\
& =&\frac{2e}{\pi}\sum_{n\geq0}\frac{(-1)^n}{(2n+1)!} \sum_{k\geq0}\frac{t^{2k}}{(2k)!} \int_0^1 e^{-x^2} x^{2n+2k+2}(1-x^2)^{n+1/2} dx.\label{Idva}
\end{eqnarray}

Now, \cite[Section 3.383(1)]{GR} implies that
\[
\int_0^1e^{-x^2} x^{2n+2k}(1-x^2)^{n + \f{1}{2}}dx=
\frac{\g(n+\f{3}{2})\g(n+k+\fs)}{2e\g(2n+k+2)} {}_1F_1(n+\f{3}{2};2n+k+2;1),
\]
so that, with $\chi=t^2/4$ the inner sum over $k$ appearing in the 
expression that emerges on the right-hand side of (\ref{Iodin}) and defines $\mathbf{I}_1$ is
\[ \frac{\sqrt{\pi}}{2e}\sum_{k\geq0}\frac{\chi^k \g(n+k+\fs)}{k! \g(k+\fs)}\sum_{r\geq0}\frac{\g(n+r+\f{3}{2})}{r! \g(2n+k+r+2)}\hspace{5cm}\]
\begin{eqnarray*}
&=&\frac{\sqrt{\pi}}{2e}\sum_{r\geq0}\frac{\g(n+r+\f{3}{2})}{r!} \sum_{k\geq0}\frac{\chi^k}{k!}\,\frac{(k+\fs)\ldots (k+n-\fs)}{\g(2n+k+r+2)}\\
&=&\frac{\sqrt{\pi}}{2e}\sum_{r\geq0}\frac{\g(n+r+\f{3}{2})}{r!}\,\chi^{1/2}\bl(\frac{d}{d\chi}\br)^n \sum_{k\geq0}\frac{\chi^{n+k-1/2}}{k! \g(2n+k+2+2)}\\
&=&\frac{\sqrt{\pi}}{2e}\sum_{r\geq0}\frac{\g(n+r+\f{3}{2})}{r!}\,\chi^{1/2} \bl(\frac{d}{d\chi}\br)^n [\chi^{-r/2-1} I_{2n+r+1}(2\sqrt{\chi})].
\end{eqnarray*}
Hence,
\[
{\bf I}_1=\sqrt{\frac{\chi}{\pi}}\sum_{n\geq0}\frac{(-1)^n}{(2n)!} \sum_{r\geq0}\frac{\g(n+r+\f{3}{2})}{r!} \bl(\frac{d}{d\chi}\br)^n [\chi^{-r/2-1} I_{2n+r+1}(2\sqrt{\chi})].
\]
A similar argument produces the representation
\[
{\bf I}_2=\sqrt{\frac{\chi}{\pi}}\sum_{n\geq0}\frac{(-1)^n}{(2n+1)!} \sum_{r\geq0}\frac{\g(n+r+\f{3}{2})}{r!} \bl(\frac{d}{d\chi}\br)^{n+1} [\chi^{-r/2-1/2} I_{2n+r+2}(2\sqrt{\chi})].
\]

This then leads to the following expressions in terms of the derivatives of specific modified Bessel functions   
for ${\bf I}_1$ and ${\bf I}_2$:
\bee		\label{konech1}
{\bf I}_1 = \frac{2t}{\sqrt{\pi}} \sum_{n\geq0}\frac{(-2)^n}{(2n)!} \sum_{r\geq0}\frac{2^r\g(n+r+\f{3}{2})}{r!} \bl(\frac{1}{t}~ \frac{d}{dt}\br)^n [t^{-r-2} I_{2n+r+1}(t)] 
\ee
and
\bee		\label{konech2}
{\bf I}_2 = \frac{2t}{\sqrt{\pi}} 
\sum_{n\geq0}\frac{(-2)^n}{(2n+1)!} \sum_{r\geq0}\frac{2^r\g(n+r+\f{3}{2})}{r!}
\bl(\frac{1}{t}  \frac{d}{dt}\br)^{n+1} [t^{-r-1} I_{2n+r+2}(t)].
\ee

To conclude, note that the validity of representation  (\ref{okonchatPRED})  is easily obtained by 
combining (\ref{inverseLaTr}), (\ref{samInteg}), (\ref{raznostINTEG}), (\ref{konech1}) and (\ref{konech2}). $\Box$

\vspace{0.6cm}

\begin{center}
\noindent {\bf 3. \  Two lemmas pertaining to a certain integrated tail which involves a Bessel function}
\end{center}
\setcounter{section}{3}
\setcounter{equation}{0}
\renewcommand{\theequation}{\arabic{section}.\arabic{equation}}

\newtheorem{lemma}{Lemma}

In this section, we establish two analytic assertions for the ``integrated tail" that involves some 
modified Bessel functions.

\begin{lemma}		\label{Mera!Levy2}
\noindent For arbitrary fixed parameter $\rho>1/2$ and argument $y \in \mathbf{R}_+^1$,
\[
2^{\rho -1} \cdot \Gamma(\rho)  \cdot \int_y^\infty \frac{e^{-x} I_{\rho-1}(x)}{x^{\rho}}  dx
\]
\begin{equation}	\label{dljaVtor1}
= 
-\gamma + y \cdot {}_3F_3\bl(\begin{array}{c} 1,1,\rho+\fs\\2,2,2\rho\end{array};-2y\br)-\log (2y)+\psi(2\rho-1)-\psi(\rho-\fs).
\end{equation}
\end{lemma}

\noindent {\bf Proof of Lemma \ref{Mera!Levy2}.}
By \cite[(10.32.2)]{DLMF}, we have the following integral representation valid for $\rho>1/2$:
\[
I_{\rho-1}(x)=\frac{(\fs x)^{\rho-1}}{\sqrt{\pi} \g(\rho-\fs)} \int_{-1}^1 (1-t^2)^{\rho-3/2} e^{xt} dt.
\]
Consider the integral
\[
{\bf J} := 2^{\rho-1} \g(\rho) \int_y^\infty \frac{e^{-x} I_{\rho-1}(x)}{x^\rho}\,dx\]
\[= \frac{\g(\rho)}{\sqrt{\pi}\g(\rho-\fs)} \int_{-1}^1(1-t^2)^{\rho-3/2} \int_y^\infty \frac{e^{-x(1-t)}}{x}\,dx\, dt,\]
where in view of \cite[(6.6.2)]{DLMF}, the inner integral equals 
\[
E_1(y(1-t)) = -\gamma-\log\,y(1-t)-\sum_{n\geq 1} \frac{(-)^n y^n(1-t)^n}{n n!}~.
\]`
Here, $E_1(z) := \int_z^{\infty} t^{-1} e^{-t} dt$ is the {\it exponential integral}. Subsequently,
\[
{\bf J} = \frac{\g(\rho)}{\sqrt{\pi}\g(\rho-\fs)} \bl\{(-\gamma-\log\,y)\int_{-1}^1 (1-t^2)^{\rho-3/2}dt\]
\[ -\int_{-1}^1 (1-t^2)^{\rho-3/2} \log (1-t)\,dt  -{\bf S} \br\}\]
\[= -\gamma -\log\,y -\frac{1}{2}\{\psi(\rho-\fs)-\psi(\rho)\}-\frac{\g(\rho) {\bf S}}{\sqrt{\pi}\g(\rho-\fs)},\]
where
\begin{eqnarray*}
{\bf S} &=&\sum_{n\geq 1} \frac{(-y)^n}{n \cdot n!} \int_{-1}^1(1-t^2)^{\rho-3/2} (1-t)^n dt\\
&=&\sum_{n\geq 1} \frac{(-y)^n}{n \cdot n!} \int_{-1}^1 (1-t)^{\rho+n-3/2} (1+t)^{\rho-3/2} dt\\
&=&\sum_{n\geq 1} \frac{(-y)^n}{n \cdot n!} 2^{n+2\rho-2} \int_0^1 w^{\rho+n-3/2} (1-w)^{\rho-3/2} dw\\
&=&2^{2\rho-2}\g(\rho-\fs)\sum_{n\geq 1} \frac{(-2y)^n}{n \cdot n!}\,\frac{\g(n+\rho-\fs)}{\g(n+2\rho-1)}\\
&=& \frac{2^{2\rho-2}\g(\rho-\fs)\g(\rho+\fs)}{\g(2\rho)} \cdot \sum_{n\geq 0}
\frac{(-2y)^{n+1} (\rho+\fs)_n}{(n+1) (2)_n (2\rho)_n}.
\end{eqnarray*}

Hence, 
\[
\frac{\g(\rho) {\bf S}}{\sqrt{\pi}\g(\rho-\fs)}=-y \sum_{n\geq0} \frac{(-2y)^n (1)_n(1)_n (\rho+\fs)_n}{(2)_n(2)_n (2\rho)_n n!}
=-y\, {}_3F_3\bl(\begin{array}{c}1,1,\rho+\fs\\2,2,2\rho\end{array};-2y\br).
\]
Using the result
\[
\psi(2\rho-1)=\fs\{\psi(\rho-\fs)-\psi(\rho)\}+\log\,2,
\]
we finally find that for $y>0$,
\[
{\bf J} = -\gamma-\log\,2y+\psi(2\rho-1)-\psi(\rho-\fs)+y\,{}_3F_3\bl(\begin{array}{c}1,1,\rho+\fs\\2,2,2\rho\end{array};-2y\br),
\]
and (\ref{dljaVtor1}) is thus proved. $\Box$

\medskip 

\begin{lemma}		\label{Poincare}
\noindent  For an arbitrary fixed $\rho>1/2$, the following Poincar\'{e} asymptotic expansion holds as $y\to+\infty$:
\bee\label{e1}
\int_y^\infty e^{-x}\,\frac{I_{\rho-1}(x)}{x^\rho}\,dx\sim \frac{y^{-\rho+1/2}}{\sqrt{2\pi}}\sum_{k=0}^\infty 
\frac{(-1)^k \Gamma(\rho+k-\fs)}{(\rho+k-\fs) \g(\rho-k-\fs)k!}\,(2y)^{-k}.
\ee 
\end{lemma}

\noindent {\bf Proof of Lemma \ref{Poincare}.}  
In order to establish (\ref{e1}), we find from  (\ref{dljaVtor1}) that
for $y \in \mathbf{R}_+^1$, 
\begin{equation}	\label{VveL}
\int_y^\infty e^{-x} \frac{I_{\rho-1}(x)}{x^\rho}\,dx = \frac{2^{1-\rho}}{\Gamma(\rho)}\bl\{-\Lambda(y)+y\cdot{}_3F_3\bl(\begin{array}{c} 1,1,\rho+\fs\\2,2,2\rho\end{array};-2y\br)\br\},
\end{equation}
where 
\[
\Lambda(y):=\gamma+\log (2y)+\psi(\rho-\fs)-\psi(2\rho-1).
\]
%
From\cite[Section 7.2.3 (12)]{PBM}, we have the integral representation
\[
{}_3F_3\bl(\begin{array}{c} 1,1,\rho+\fs\\2,2,2\rho\end{array};-2y\br)=\frac{\g(2\rho)}{\g(\rho+\fs)}\,\frac{1}{2\pi i}\int_L \frac{\g(-s) \g(s+\rho+\fs)}{(s+1)^2 \g(s+2\rho)}\,(2y)^s ds,
\]
where $L$ is a path parallel to the imaginary $s$-axis that separates the poles of $\g(-s)$ from the double pole at $s = -1$ and the simple (when $\rho> 1/2$) poles at $s= -\rho-\fs-k$, where 
$k=0, 1, 2 \ldots $. Noting that the residue at the double pole $s = -1$ is given by $\Lambda(y)/y$, we find upon  displacement of the integration path to the left in the usual manner that as $y\to+\infty$,
\[
y\, {}_3F_3\bl(\begin{array}{c} 1,1,\rho+\fs\\2,2,2\rho\end{array};-2y\br)\sim 
\Lambda(y)+\frac{y \g(2\rho)}{\g(\rho+\fs)}\hspace{5cm}\]
\[\hspace{3cm} \times \sum_{k=0}^\infty \frac{(-1)^k \g(\rho+k-\fs)}{(\rho+k-\fs) \g(\rho-k-\fs)}\,(2y)^{-\rho-k-1/2}.\]
The representation (\ref{e1}) then follows upon use of the duplication formula for the gamma function. $\Box$


\vspace{0.6cm}

\begin{center}
\noindent {\bf 4. \  Two modified Bessel function inequalities}
\end{center}
\setcounter{section}{4}
\setcounter{equation}{0}
\renewcommand{\theequation}{\arabic{section}.\arabic{equation}}
In this section, we establish the following two (possibly) unknown inequalities involving the modified Bessel function $I_\nu(x)$ whose consideration was motivated by 
the needs of probability theory, but which are also important for the theory of special functions in their own right (see formulas (\ref{a1}) and (\ref{a2})).

\begin{theorem}		\label{NovajaVazhnaja}
Fix an arbitrary real $\nu> -1/2$ and assume that the argument $x \in \mathbf{R}_+^1$. Then the following two 
inequalities hold: 
\begin{equation}		\label{a1}
\frac{1}{x} I_\nu(x) I_{\nu+1}(x)>I_{\nu+1}(x)^2-I_\nu(x) I_{\nu+2}(x);
\end{equation}
%
%
\begin{equation}		\label{a2}
I_\nu(x)>\bl(1+\frac{2\nu+1}{2x}\br) I_{\nu+1}(x).
\end{equation}
\end{theorem}

The second inequality (\ref{a2}) implies that for $x>0$ and $\nu > -1/2$,
\begin{equation}		\label{a3}
I_\nu(x) > I_{\nu+1}(x),
\end{equation}
which was established by Jones  \cite{Jones}. 

\medskip

\noindent{\bf Proof of Theorem  \ref{NovajaVazhnaja}.}\\
(i) In order to prove (\ref{a1}), we denote $f(x) :=e^{-x}g(x)$, where 
\begin{equation}		\label{vspom}
g(x) := x^{-\nu} I_\nu(x). 
\end{equation}
Then 
\begin{equation}		\label{UPROSTa1DoublePrime}
f^{\prime}(x)=e^{-x}(g^{\prime}(x) - g(x)) = e^{-x} x^{-\nu} (I_{\nu + 1}(x) - I_\nu(x))<0
\end{equation}
by virtue of (\ref{a3}). 

\medskip

 Also,
\[
f^{\prime \prime}(x)=e^{-x}(g(x)-2g^{\prime}(x) +g^{\prime \prime}(x)),
\]
and so the convexity condition $f^{\prime \prime}(x) f(x) > (f^{\prime}(x))^2$ becomes
\[ e^{-2x}g(x)(g(x)-2 g^{\prime}(x) + g^{\prime \prime}(x)) >e^{-2x} (g^{\prime}(x)-g(x))^2\]
\[ \qquad \qquad = e^{-2x}(g(x)^2+g^{\prime}(x)^2 - 2 g(x) g^{\prime}(x)),\]
to yield
\begin{equation}	\label{a2REPEAT}
g(x) \cdot g^{\prime \prime}(x) > g^{\prime}(x)^2.
\end{equation}
Now  
\[
g^{\prime}(x) =x^{-\nu}I_{\nu+1}(x),~~~~~\mbox{and}~~~~~g^{\prime \prime}(x)=x^{-\nu-1}I_{\nu+1}(x)+x^{-\nu}I_{\nu+2}(x),
\]
so that the convexity condition (\ref{a2REPEAT}) reduces to the inequality (\ref{a1}).
%

From \cite[(10.32.15)]{DLMF}, we have the integral representation
\[
I_\mu(x) I_\nu(x)=\frac{2}{\pi}\int_0^{\pi/2} I_{\mu+\nu}(w) \cos (\mu-\nu)\theta\,d\theta
\]
for $\mu+\nu > -1$, where $w := 2x\cos \theta$. This enables us to express the products of Bessel functions in inequality (\ref{a1}) in 
terms of a single Bessel function.
Then (\ref{a1}) becomes $J>0$, where
\[
J:=\frac{2}{\pi} \int_0^{\pi/2} \{I_{2\nu+1}(w) \cos \theta- x I_{2\nu+2}(w)(1-\cos 2\theta)\}\,d\theta
\]
\[
\hspace{0.4cm}=\frac{2}{\pi}
\int_0^{\pi/2}\{I_{2\nu+1}(w) \cos \theta-2x I_{2\nu+2}(w) \sin^2\theta\} d\theta.
\]
An integration by parts applied to the first integral then yields
\[
J=\frac{2}{\pi}I_{2\nu+1}(0)+\frac{4x}{\pi}\int_0^{\pi/2} \{I'_{2\nu+1}(w)-I_{2\nu+2}(w)\} \sin^2\theta\,d\theta.
\]
But
\[
I^{\prime}_{2\nu+1}(w)=I_{2\nu+2}(w)+\frac{2\nu+1}{w} I_{2\nu+1}(w)
\]
and so
\[
J=\frac{2}{\pi}I_{2\nu+1}(0)+\frac{2(2\nu+1)}{\pi}\int_0^{\pi/2} I_{2\nu+1}(w)\,\frac{\sin^2\theta}{\cos \theta}\,d\theta.
\]
The integral behaves near $\theta= \pi/2$ like $\cos^{2\nu+1}\theta/\cos \theta$ and so converges provided $\nu> -1/2$ 
and in this case $I_{2\nu+1}(0)=0$. Since the integrand is non-negative in $\theta\in[0,\pi/2]$ it follows that $J>0$, 
which concludes the proof of the validity of inequality (\ref{a1}).\\
%
\medskip

(ii) In order to prove inequality (\ref{a2}), note that with $f(x)$ defined as above, we have
\[ 
f^{\prime \prime}(x)=e^{-x}(g(x)-2g^{\prime}(x)+g^{\prime \prime}(x))=x^{-\nu}e^{-x}\{I_\nu(x)+I_{\nu+2}(x)\]
\[ + x^{-1} I_{\nu+1}(x)-2I_{\nu+1}(x)\}
= 2x^{-\nu}e^{-x}\bl\{I_\nu(x)-\bl(1+\frac{2\nu+1}{2x}\br)\,I_{\nu+1}(x)\br\}\]
upon use of the fact that
\[
I_\nu(x)+I_{\nu+2}(x)=2I_\nu(x)-\frac{2\nu+2}{x}~ I_{\nu+1}(x).
\]

From \cite[(10.32.2)]{DLMF} we have the integral representation
\[
I_\nu(x)=\frac{(\fs x)^\nu}{\sqrt{\pi} \g(\nu+\fs)} \int_0^\pi e^{x\cos \theta} \sin^{2\nu} \theta \,d\theta \qquad (\nu > -1/2)
\]
so that
\[
J^{\prime}:=\sqrt{\pi} \cdot \g(\nu+1/2) (x/2)^{-\nu}\bl\{I_\nu(x)-\bl(1+\frac{2\nu+1}{2x}\br)\,I_{\nu+1}(x)\br\}\hspace{4cm}
\]
\begin{eqnarray*}&=& \int_0^\pi e^{x\cos \theta} \sin^{2\nu} \theta \bl\{1-\bl(1+\frac{2\nu+1}{2x}\br)\,\frac{x\sin^2\theta}{2\nu+1}\br\}d\theta\\
&=&\int_0^\pi e^{x\cos \theta} \sin^{2\nu} \theta \cdot (1-\fs \sin^2\theta)\,d\theta-\frac{x}{2\nu+1}\int_0^\pi e^{x\cos \theta} \sin^{2\nu+2}\theta\,d\theta\br\}.
\end{eqnarray*}
An integration by parts applied to the second integral yields
\[
\frac{x}{2\nu+1}\int_0^\pi e^{x\cos \theta} \sin^{2\nu+2}\theta\,d\theta=\int_0^\pi e^{x\cos \theta} \sin^{2\nu} \theta\,\cos \theta\,d\theta,
\]
where the integrated term vanishes provided $\nu > -1/2$. Hence,
\begin{eqnarray*}
J^{\prime}&=& \int_0^\pi e^{x \cos \theta} \sin^{2\nu} \theta \cdot \{1-\fs \sin^2 \theta-\cos \theta \} \,d\theta\\
&=&\int_0^\pi e^{x\cos \theta} \sin^{2\nu} \theta \cdot \sin^4 (\theta/2)\,d\theta>0
\end{eqnarray*}}
thereby establishing the validity of (\ref{a2}). $\Box$ 

\bigskip

\noindent\textbf{Acknowledgements}\\

\noindent VVV is grateful to Huaxiong Huang, Neil Madras and Tom Salisbury for their
warm hospitality and thanks the Fields Institute, University of Toronto and York University for continued 
support.
 
\vspace{0.6cm}

\end{document}